\documentclass{article}
\usepackage{amssymb}
\usepackage{amsmath}

\newtheorem{theo}{Theorem}

\newtheorem{lemm}{Lemma}

\newtheorem{rema}{Remark}

\newcommand{\cqfd}
{%
\mbox{}%
\nolinebreak%
\hfill%
\rule{2mm}{2mm}%
\medbreak%
\par%
}
\newfont{\gothic}{eufb10}

\title{Torsion cohomology classes and algebraic cycles on complex projective manifolds}
\author{C. Soul\'e, C. Voisin}
\date{ \ }

\begin{document}

\maketitle

Let $X$ be a smooth complex projective manifold and $H^n (X,{\mathbb Z})$ its singular
 cohomology group of degree $n$,
with integral coefficients. Given a torsion class $\alpha \in H^{2k} (X,{\mathbb Z})$, can
we say that this class $\alpha$ is algebraic?

This is true when $k=1$, and, apparently, Hodge thought that this would always be the case \cite{H}. However,
 Atiyah and Hirzebruch found counterexamples to Hodge's assertion \cite{AH}. This
  is why the Hodge conjecture is now formulated for rational cohomology classes only. Recently,
  Totaro gave a new interpretation of the Atiyah-Hirzebruch counterexamples, in terms
  of the complex cobordism ring
   of $X$ \cite{totaro}.

When looking at these examples, we noticed that the order of the non-algebraic torsion class $\alpha$
can be divisible only by primes at most equal to the complex dimension of $X$ (Th. 1). On the other hand,
a construction of Koll\'ar \cite{kollar} provides examples of manifolds $X$ as above with
 a class $\alpha \in H^4 (X,{\mathbb Z})$ which is not algebraic, while a non-zero high multiple of $\alpha$ is
 algebraic. Inspired by this construction, we give examples of a non-algebraic $p$-torsion class
 in $H^6 (X,{\mathbb Z})$,
 with $\dim_{\mathbb C} (X) = 5$ and with $p$ any prime bigger than 3 (Th. 3).
  In particular, these classes could not be studied
 by the (topological) methods of \cite{AH} and \cite{totaro}.

In \cite{totaro}, Totaro gives examples of non trivial torsion classes in the
 Griffiths group of homologically trivial algebraic cycles
modulo those algebraically equivalent to zero. Such examples were
also constructed by Schoen \cite{schoen},
 but Totaro's construction provides non trivial torsion
 cycles annihilated
 by the Abel-Jacobi map, while Schoen uses the Abel-Jacobi invariant to conclude
 as in  Griffiths \cite{griffiths}  that his cycles are
 not algebraically equivalent to $0$. Note that for codimension $2$ cycles, the Abel-Jacobi map
 is known to be injective on torsion cycles homologous to $0$
 \cite{colliot}.
 Theorem 4, again inspired by Koll\'ar's argument, gives new examples of such algebraic cycles,
which furthermore cannot be detected by Totaro's method (nor indeed by any locally constant invariant).
More generally, we get non trivial algebraic cycles
in any level of the Hiroshi Saito  filtration on Chow groups \cite{hsaito} (Theorem  \ref{th5}).

On the positive side, Bloch made the beautiful remark \cite{B} that a conjecture of his
and Kato on the Milnor $K$-theory of fields \cite{BK} implies that any torsion class
in $H^n (X,{\mathbb Z})$, $n > 0$, is supported in codimension one. We note in Theorem \ref{th7}
that it implies also that the image of this class by the Atiyah-Hirzebruch differentials
are supported in codimension two.

Sections 1 and 5 (resp. 2, 3, 4) are due to the first (resp. second) author.
This work started at the Hodge's Centennial Conference in Edinburgh, 2003.
 Both authors thank the organizers of that meeting for their invitation. We also thank
P. Guillot and B. Totaro for useful comments, and
J. Koll\'ar for allowing us to reproduce his arguments in section 2.

\section{On the counterexamples of Atiyah and Hirzebruch}\label{sec1}

Let $X$ be a smooth projective complex manifold, $H^n (X,{\mathbb Z})$ its codimension $n$ singular
 cohomology with integral coefficients,
and ${\rm CH}^k (X)$ the group of codimension $k$ algebraic cycles on $X$ modulo rational equivalence.
Any codimension $k$ cycle $Z$ on $X$
defines a cohomology class $[Z] \in H^{2k} (X,{\mathbb Z})$, the image of which in $H^{2k} (X,{\mathbb C})$ has Hodge type $(k,k)$.
 Hodge asked whether, conversely, any class $\alpha \in H^{2k} (X,{\mathbb Z})$ with
  image $\alpha_{\mathbb C} \in H^{2k} (X,{\mathbb C})$ of
 type $(k,k)$ is of the form $\alpha = [Z]$ for some cycle $Z$ on $X$ (\cite{H}, end of \S~2). In particular,
 when $\alpha$ is a {\it torsion class}, {\it i.e.} $\alpha_{\mathbb C} = 0$, is it algebraic?
 This is true when $k=1$, but Atiyah and Hirzebruch found counterexamples when $k \geq 2$ \cite {AH}.
 In these examples, to prove that $\alpha \in H^{2k} (X,{\mathbb Z})$ is not algebraic,
  they use the following criterion. Consider the spectral sequence of generalized cohomology, with
$$
E_2^{ st} (X) = \left\{ \begin{matrix} H^s (X,{\mathbb Z}) &\hbox{if $t$ is even} \\ 0 &\hbox{if $t$ is odd,} \end{matrix} \right.
$$
which converges to the topological $K$-groups $K_{\rm top}^{s + t} (X)$. If a class $\alpha \in H^{2k} (X,{\mathbb Z})$ is algebraic,
its image by every differential $d^r$, $r \geq 2$, in that spectral sequence must vanish (\cite{AH}, Th. 6.1).

Recently, Totaro revisited these examples of Atiyah and Hirzebruch. Let $MU^* (X)$ be
the complex cobordism graded ring of $X$, and $MU^* (X) \underset{MU^*}{\otimes} {\mathbb Z}$ its tensor product with ${\mathbb Z}$ over $MU^* = MU^*(\hbox{point})$
 (which maps to ${\mathbb Z} = H^0 (\hbox{point}, {\mathbb Z})$). Totaro proved that the cycle map
$$
{\rm CH}^k (X) \to H^{2k} (X,{\mathbb Z})
$$
is the composite of two maps
\begin{eqnarray}\label{mapnom}
{\rm CH}^k (X) \to (MU^* (X) \underset{MU^*}{\otimes} {\mathbb Z})^{2k} \to H^{2k} (X,{\mathbb Z}) \, .
\end{eqnarray}
Therefore, a torsion class which is not
in the image of $(MU^* (X) \underset{MU^*}{\otimes} {\mathbb Z})^{2k}$ cannot be algebraic.

Our first remark is that, when one uses the Atiyah-Hirzebruch
 criterion or Totaro's factorization of the cycle map, the examples of non-algebraic torsion classes
  one gets must be of small order compared to the dimension of $X$.

\medskip

\begin{theo} {\it Let $p $ be a prime integer and $\alpha \in H^{2k} (X,{\mathbb Z})$ a cohomology class
such that $p \alpha = 0$. Assume that $p$ is bigger than the complex dimension of $X$. Then,
for every $r \geq 2$, $d^r (\alpha) = 0$. Furthermore, $\alpha$ lies in the image of
\begin{eqnarray}\label{formulephik}
\phi^k:(MU^* (X) \underset{MU^*}{\otimes} {\mathbb Z})^{2k} \to H^{2k} (X,{\mathbb Z}) \, .
\end{eqnarray}
}
\end{theo}

\noindent {\bf Proof.} The first assertion is rather standard. For any $q > 1$, the
$q$-th Adams operation $\psi^q$ acts upon the Atiyah-Hirzebruch spectral sequence $E_r^{ st} (X)$.
When $t = -2i$, its action on
$E_2^{ st} (X)$ is the multiplication by $q^i$. Since $\psi^q$ commutes with $d^r$, we get
$$
(q^{k+r_0} - q^k) \, d^r (\alpha) = 0
$$
when $r = 2r_0 + 1$ ($d^r$ is zero when $r$ is even). For the prime number $p$ to divide $q^{k+r_0} - q^k$
 for every $q > 1$,
 it is necessary that $p-1$ divides $r_0$. On the other hand, since $p\alpha = 0$, we must have $k > 0$
 and $d^r (\alpha)$ is not of top degree. Therefore $r \leq 2 \dim_{\mathbb C} (X) - 2$,
 hence $r_0 < \dim_{\mathbb C} (X)$. Since $p > \dim_{\mathbb C} (X)$, we get a contradiction.

Now consider the edge homomorphism
$$
e^k : K_{\rm top}^0 (X) \to H^{2k} (X,{\mathbb Z})
$$
in the Atiyah-Hirzebruch spectral sequence. Since $d^r (\alpha) = 0$ for every $r \geq 2$, we have
\begin{equation}
\label{eq31}
\alpha = e^k (\xi)
\end{equation}
for some virtual bundle $\xi \in K_{\rm top}^0 (X)$. We want to show that $\alpha$ lies in the image
of the map
$$\phi^k : (MU^* (X) \underset{MU^*}{\otimes} {\mathbb Z})^{2k} \to H^{2k} (X,{\mathbb Z}) \, .$$

First notice that, for every $i \geq 0$, the Chern class $c_i (\xi) \in H^{2i} (X,{\mathbb Z})$ is in
 the image of $\phi^i$. Indeed we have
\begin{equation}
\label{eq32}
c_i (\xi) = \phi^i ({\rm cf}_{\alpha} (\xi) \otimes 1) \, ,
\end{equation}
where ${\rm cf}_{\alpha} (\xi) \in MU^{2i} (X)$ is the $i$-th Conner-Floyd class of $\xi$,
 with $\alpha = (i,0,0,\ldots)$ (\cite{A}, Th. 4.1).
 To check (\ref{eq32}), by the splitting principle, we can assume that $i=1$, in which case it follows from the definitions (loc. cit.).

Since $\phi^*$ is a ring homomorphism, any polynomial in the Chern classes $c_i (\xi)$ lies
in its image. Let $N_k (c_1 , \ldots , c_k)$ be the $k$-th Newton polynomial. We claim that
\begin{equation}
\label{eq33}
k! \, e^k (\xi) = N_k (c_1(\xi) , \ldots , c_k (\xi)) \, .
\end{equation}
To check (\ref{eq33}), since $\xi$ is the pull-back of a vector bundle on a Grassmannian, we can assume
that $H^* (X,{\mathbb Z})$ is torsion free. Then it is enough to prove that the image of
$e^k (\xi)$ in the rational cohomology $H^{2k} (X,{\mathbb Q})$ is
$$
{\rm ch}_k (\xi) = N_k (c_1 (\xi) ,  \ldots , c_k (\xi)) / k! \, .
$$
But this identity follows from the fact that the Chern character
$$
{\rm ch} : K_{\rm top}^* (-) \to H^* (- , {\mathbb Q})
$$
is a morphism of extraordinary cohomology theories (\cite{AH}, \S~2).

Since $p\alpha = 0$ with $p > \dim_{\mathbb C} (X) \geq k$,
and since $k! \, e^k (\xi)$ lies in the image
of $\phi^k$ by (\ref{eq32}), we deduce from (\ref{eq31}) and (\ref{eq33}) that $\alpha$ is also
in the image of $\phi^k$.

\section{An argument due to  Koll\'ar\label{sec1}}
We start this section by describing a method due to Koll\'ar \cite{kollar},
which produces examples
of smooth projective complex varieties $X$, together with an even degree integral cohomology
class $\alpha$, which is not algebraic, that is, which is not the cohomology class
of an algebraic cycle  of $X$, while  a non-zero multiple of $\alpha$ is algebraic.
This is another sort of counterexample to the Hodge conjecture over the integers, since
the class $\alpha$ is of course a Hodge class, the other known examples being
that of torsion classes \cite{atiyah} that we shall revisit in section \ref{sec2}.

The examples are as follows : consider a smooth hypersurface
$X\subset \mathbb{P}^{n+1}$ of degree $D$. For $l<n$ the Lefschetz theorem on hyperplane
sections says that the restriction map
$$H^l(\mathbb{P}^{n+1},\mathbb{Z})\rightarrow H^l(X,\mathbb{Z})$$
is an isomorphism.
Since the left-hand side is isomorphic to
$\mathbb{Z}H^k$ for $l=2k<n$, where $H$ is the
cohomology class of a hyperplane, and $0$ otherwise, we conclude by Poincar\'e duality on
$X$ that for $2k>n$, we have $H^{2k}(X,\mathbb{Z})=\mathbb{Z}\alpha$, where
$\alpha$ is determined by the condition $<\alpha,h^{n-k}>=1$,
with the notation
$h=H_{\mid X}=c_1(\mathcal{O}_X(1))$. Note that the class $D\alpha$ is equal to $h^{k}$,
(both have intersection number $D$ with  $h^{n-k}$),  hence
is algebraic.

In the sequel, we consider for simplicity the case where $n=3,\,k=2$. Then
$D\alpha$ is the class of a plane section of $X$.
\begin{theo}(Koll\'ar, \cite{kollar}) Assume that for some  integer
$p$ coprime to $6$, $p^3$ divides $D$. Then for general $X$, any curve $C\subset X$ has
degree divisible by $p$. Hence the class $\alpha$ is not algebraic.
\end{theo}
Recall that ``general'' means that the defining equation for $X$ has to be chosen
away from a specified  union of countably many Zariski closed proper subsets of the parameter
space.

\vspace{0,5cm}

\noindent {\bf Proof.} Let $D=p^3s$, and let $Y\subset \mathbb{P}^4$ be a degree $s$
smooth hypersurface.
Let $\phi_0,\ldots,\phi_4$ be sections of $\mathcal{O}_{\mathbb{P}^4}(p)$
 without common zeroes. They provide a map
 $$\phi:Y\rightarrow \mathbb{P}^4,$$
 which for a generic choice of the $\phi_i$'s satisfies the following properties :
 \begin{enumerate}
 \item $\phi$ is generically of degree $1$ onto its image, which is a hypersurface
 $X_0\subset \mathbb{P}^4$ of degree $p^3s=D$.
 \item \label{item2}
 $\phi$ is two-to-one generically over a surface in $X_0$, three-to-one
 generically over a curve  in $X_0$,  at most
 finitely many points of $X_0$ have more than $3$ preimages, and no point has more than
 $4$ preimages.

 \end{enumerate}
 Let $\mathbb{P}^N$ be the projective space of all polynomials of degree $D$ on $\mathbb{P}^4$,
 and let $\mathcal{X}\rightarrow\mathbb{P}^N$ be the universal hypersurface.
 Introduce the relative Hilbert schemes (cf \cite{kollarbook})
 $$\mathcal{H}_\nu\rightarrow \mathbb{P}^N,$$
 parameterizing pairs $\{(Z,X),Z\subset X\}$, where $Z$ is a $1$-dimensional subscheme
 with Hilbert polynomial $\nu$. The Hilbert polynomials
 $\nu$ encode the degree
 and genus of the considered subschemes, hence there are only countably many of them.
 The important points are the following :
 \begin{itemize} \item The morphism $\rho_\nu:\mathcal{H}_\nu\rightarrow \mathbb{P}^N$
 is projective.
 \item There exists a universal subscheme
 $$\mathcal{Z}_\nu\subset \mathcal{H}_\nu\times_{{\mathbb P}^N}\mathcal{X}$$
 which is flat over $\mathcal{H}_\nu$.
 \end{itemize}
 Let $U$ be the set
 $$\mathbb{P}^N\setminus\bigcup_{\nu\in I}\rho_{\nu}(\mathcal{H}_\nu),$$
 where the set $I$ is the set of Hilbert polynomials $\nu$ for which
 the map $\rho_{\nu}$ is not dominating. Let now
 $X\subset \mathbb{P}^4$ be a smooth hypersurface
 which is parameterized by a point $x\in U$ (so $X$ is general). Let $C\subset X$ be a curve.
 The reduced structure on $C$ makes $C$ into a subscheme of $X$, which is parameterized
 by a point $c_x\in \mathcal{H}_\nu$ over $x$, for some $\nu$. By definition of
 $U$, since $x=\rho_{\nu}(c_x)$, the map $\rho_{\nu}$ has to be dominating, hence surjective.
  Hence it follows that there is a point
  $c_0\in \mathcal{H}_\nu$ over the point $x_0$ parameterizing the hypersurface $X_0$.
  The fiber $Z_0$ of the universal subscheme $\mathcal{H}_\nu$ over $c_0$
  provides a subscheme $Z_0\subset X_0$, which by flatness has the same degree
  as $C$. Let $z_0$ be the associated cycle of $X_0$. Recall the normalization
  map
  $$\phi:Y\rightarrow X_0.$$
  By property \ref{item2} above, there exists a $1$-cycle
  $\tilde{z}_0$ in $Y$ such that $\phi_*(\tilde{z}_0)=6z_0$. It follows that
  $$6deg\,z_0= deg\,\phi_*(\tilde{z}_0).$$
  On the other hand, the right-hand side is equal to the degree of the line bundle
  $\phi^*\mathcal{O}_{X_0}(1)$ computed on the cycle $\tilde{z}_0$.
  Since $\phi^*\mathcal{O}_{X_0}(1)$ is equal to
  $\mathcal{O}_Y(p)$ (recall that $Y\subset \mathbb{P}^4$ was a hypersurface of degree
$s$), it follows that this degree is divisible by $p$.
Hence we found that $6deg\,C$ is divisible by $p$, and since $p$ is coprime to
$6$, it follows that $deg\,C$ is also divisible by $p$.

\cqfd
\begin{rema}\label{remarque}{\rm In contrast, one can show that there exists a countable union
of proper algebraic subsets, which is dense in the parameter space $\mathbb{P}^N$,
parameterizing hypersurfaces $X$ for which the class $\alpha$ is algebraic.
It suffices for this to prove that the set of  surfaces of degree $D$ carrying an
algebraic class
$\lambda\in H^{2}(S,\mathbb{Z})\cap H^{1,1}(S)$ satisfying the property
that $<\lambda,c_1(\mathcal{O}_S(1))>$ is coprime to $D$, is dense in the space of all
surfaces of degree $D$ in $\mathbb{P}^3$. Indeed, for any $X$ containing such a surface,
the class $\alpha$ is algebraic on $X$.

Now this fact follows from the density criterion for the Noether-Lefschetz locus explained
in \cite{voisin} 5.3.4, and from the fact that rational classes
 $\lambda\in H^2(S,\mathbb{Q})$ such that a multiple $b\lambda$ is integral,
and satisfies $<b\lambda,c_1(\mathcal{O}_S(1))>=a $ with $a$  coprime to $D$, are dense in
$H^2(S,\mathbb{Q})$.

}
\end{rema}

 To conclude this section, we note that Koll\'ar's construction
works only for high degree hypersurfaces, and indeed hypersurfaces
$X$ of degree $\leq 2n-1$ in $\mathbb{P}^{n+1}$ contain lines, whose
class is equal to the positive generator of
$H^{2n-2}(X,\mathbb{Z})$. So the following question might still have
a positive answer :

{\it Let $X$ be a Fano variety
or more generally a rationally
connected variety of dimension $n$. Is the Hodge conjecture true for the integral
cohomology classes
of degree $2n-2$ on $X$?}

Note that all such classes are of type $(n-1,n-1)$, since $H^2(X,\mathcal{O}_X)=0
=H^{n-2}(X,K_X)^*$. Since the rational Hodge conjecture
is known to be true for degree $2n-2$ classes,
 it is always true that a multiple of such a class is
algebraic. Note also that the question has a negative answer
for classes of degree $2n-2k,\,n-3\geq k\geq2$, at least in the rationally
connected case. Indeed, it suffices to start with
one of Koll\'ar's examples $X\subset \mathbb{P}^{n+1}\subset \mathbb{P}^{n+l},\,l\geq2$, and
to blow-up $X$ in $ \mathbb{P}^{n+l}$. Since $X$ has a degree $2n-2$
integral class which is of Hodge type $(n-1,n-1)$, but is not algebraic,
the resulting variety $Y$, which is $n+l$-dimensional and rationally connected, has
degree $2n+2s,\,l-2\geq s\geq 0$
integral classes which are of Hodge type $(n+s,n+s)$, but are not algebraic.

One reason to ask this question is the fact that there is the following criterion
for rationality :
\begin{lemm} Let $X$ be a variety which is birationally equivalent
to $\mathbb{P}^n$. Then any integral class of degree $2n-2$ on $X$ is algebraic.
Furthermore, the Hodge conjecture is true for degree $4$ integral Hodge classes
on $X$.
\end{lemm}
\noindent {\bf Proof.}
In both cases, there is a variety $Y$ which admits a
morphism of degree $1$, $\phi:Y\rightarrow X$, and is obtained
from $\mathbb{P}^n$ by a sequence of blow-ups along smooth centers.
Using the $\phi_*$ map, one concludes that if the statement is true for $Y$,
 it is true for $X$. Since the statement is true for $\mathbb{P}^n$, it suffices then
 to show that if the statement is true for a smooth projective variety $Z$, it is true
 for the blowing-up $Z_W$ of $Z$ along a smooth center $W$.
 But the supplementary classes of degree $2n-2$  on $Z_W$ are generated by
 classes of curves contracted by the blown-down map, hence they are algebraic.
 The supplementary integral  Hodge classes of degree $4$ on $Z_W$ come from
integral  Hodge classes of degree $2$ or $0$ on $W$, hence they are also algebraic.
\cqfd
\section{\label{sec2}Application to torsion classes}
As  said in the first section,
examples of even degree torsion cohomology classes
which are not algebraic were first found by
Atiyah and Hirzebruch \cite{atiyah}. They exhibited topological obstructions for a
torsion class to be the cohomology class of an algebraic cycle. These obstructions were
reinterpreted by Totaro \cite{totaro},
who stated  the following criterion for algebraicity :

{\it  (*)  For a  degree $2k$  class  to be algebraic,
it has  to be in the image of the map $\phi^k$ of  (\ref{formulephik}).}

This leads to the construction of even degree torsion classes
which are not algebraic. On the other hand, as shown by
Theorem 1,
 these examples must be of large dimension,
and  the criterion above cannot be applied to provide for any prime $p$,
examples of $p$-torsion classes of a given degree on projective manifolds of given dimension,
which are not algebraic.

In this section  we apply Koll\'ar's argument to construct,
for any prime $p\geq5$, examples of
$p$-torsion
cohomology classes of degree $6$ on smooth  projective varieties
$X$ of dimension
$5$, which are not algebraic.

Furthermore, the fact that these classes are not algebraic cannot
be detected by topological arguments. Indeed, any obstruction to algebraicity
which is locally constant on the parameter space of $X$
must vanish on these classes, as we show that they become algebraic on a
dense subset of the  parameter space  of $X$. In particular,
these classes are in the image of
the map $\phi^k$, which shows that criterion (*) is not sufficient.

The construction is as follows : let $p\geq5$ be a prime integer. Let $S$ be a surface
which admits a copy of $\mathbb{Z}/p\mathbb{Z}$  as a direct summand in $H^2(S,\mathbb{Z})$.
Such a surface can be constructed by a Godeaux type construction : namely
one can take for $S$ the quotient of a degree $p$
smooth surface $\Sigma$ in $\mathbb{P}^3$, defined by an equation invariant under
$g$, where
$g$ acts on homogeneous coordinates
by
$$ g^*X_i=\zeta^i X_i,\,i=0,\ldots,3.$$
Here $\zeta$ is a $p$-th root of unity, so that $g$ has order $p$.
As $\Sigma$ is simply connected, one sees easily that the torsion
of $H^2(S,\mathbb{Z})$ is isomorphic to $\mathbb{Z}/p\mathbb{Z}$, generated by
$c:=c_1(\mathcal{L})$, where $\mathcal{L}$ is any of the $p$-torsion line bundles on
$S$ corresponding to a non trivial character of
$\mathbb{Z}/p\mathbb{Z}$. Note that the class $c$ is not divisible by $p$, since
the torsion of $H^2(S,\mathbb{Z})$ is isomorphic to $\mathbb{Z}/p\mathbb{Z}$.

We now consider a hypersurface
 $X\subset \mathbb{P}^4$ of degree $p^3$. Recall from the previous
section
that $H^4(X,\mathbb{Z})$ is generated by $\alpha$, where
$\alpha$ satisfies the condition
$$<\alpha,h>=1.$$
The class $\gamma:=pr_1^*c\cup pr_2^*\alpha$ is a degree $6$ cohomology class
on $S\times X$, which is of $p$-torsion. Observe that, according to remark \ref{remarque},
this class is algebraic
for a dense set of parameters for $X$. Indeed, since
$c=c_1(\mathcal{L})$ is algebraic, that is $c=[D]$ for some divisor $D$ of $S$,
 once $\alpha$ is algebraic, say
$\alpha=[Z],\,Z\in \mathcal{Z}^2(X)$, we also have
$$\gamma=[pr_1^*D\cap pr_2^* Z].$$
We now have :
\begin{theo} For general $X$, the class $\gamma$ is not algebraic.
More precisely, for $X$ general in modulus, and for any surface
$T\subset S\times X$, the K\"unneth component $[T]^{2,4}$ of $[T]$
which lies in $H^2(S,\mathbb{Z})\otimes H^4(X,\mathbb{Z})$ is of the form
$t\otimes  \alpha$, where $t$ is divisible by $p$.
\end{theo}

(Note that the K\"unneth decomposition is well defined for
$S\times X$, because the cohomology of $X$ has no torsion.)

\vspace{0.5cm}

\noindent {\bf Proof.} For $X,\,T$ as in the theorem, let the
K\"unneth component $[T]^{2,4}$ of $[T]$ be of the form
$t\otimes\alpha$. If $j$ is the natural inclusion map
of $X$ into $\mathbb{P}^4$, it follows that the K\"unneth component of
$[(Id,j)(T)]$ which is of type $(2,6)$ is equal to
$t\otimes j_*\alpha$. Since the class $j_*\alpha$ generates
$H^6(\mathbb{P}^4,\mathbb{Z})$, the statement is equivalent to the fact that
$[(Id,j)(T)]^{2,6}$ is divisible by $p$.

We now apply the argument of section \ref{sec1}, except that instead of considering
the relative Hilbert schemes of $1$-dimensional subschemes in hypersurfaces $X$, we
consider the
relative Hilbert schemes of $2$-dimensional subschemes in the products $S\times X$.

This provides us with a set $U$ of parameters, which is the complementary set
of a countable union of proper algebraic subsets in $\mathbb{P}^N$, and has the property
that any surface $T\subset S\times X$ admits a flat specialization
$T_0\subset S\times X_0$, for any  specialization of $X$.

We will take  for $X_0$ the image of a generic morphism
$$\phi:\mathbb{P}^3\rightarrow \mathbb{P}^4$$
given by a base-point free linear system of polynomials degree $p$.

Now, by property \ref{item2} (see section \ref{sec1}) of the morphism $\phi$,
there exists a $2$-cycle $\tilde{t}_0$ in
$S\times \mathbb{P}^3$, such that, denoting by $t_0$ the cycle associated to
the subscheme
$T_0$, we have the following equality of cycles in $S\times X_0$ :
$$12t_0=(Id,\phi_0)_*\tilde{t}_0,$$
where $\phi_0$ is the map $\phi$ viewed as a map from
$\mathbb{P}^3$ to $X_0$.
It follows that we also have  the equality of $2$-cycles in $S\times\mathbb{P}^4$:
$$12(Id,j_0)_*t_0=(Id,\phi)_*\tilde{t}_0,$$
where $j_0$ is the inclusion of $X_0$
in $\mathbb{P}^4$.
This equality translates into an equality between
cycle classes :
$$12[(Id,j_0)_*(t_0)]=[(Id,\phi)_*\tilde{t}_0],$$

and between their K\"unneth components of type $(2,6)$ :
\begin{eqnarray}\label{formule}12[(Id,j_0)_*t_0]^{2,6}=[(Id,\phi)_*\tilde{t}_0]^{2,6}.
\end{eqnarray}
Now note that since $T_0$ is a flat specialization of $T$,
we have $[(I,j_0)(T_0)]=
[(Id,j)(T)]$ and hence
$$[(Id,j_0)_*t_0]^{2,6}=[(Id,j)(T)]^{2,6}=t\otimes j_*\alpha.$$
So the left-hand side in (\ref{formule}) is equal to $12t\otimes j_*\alpha$.
On the other hand, the right-hand side in (\ref{formule})
is equal to
$(Id,\phi)_*([\tilde{t}_0]^{2,4})$, where $[\tilde{t}_0]^{2,4}$ is the K\"unneth
component of type
$(2,4)$ of the class
$[\tilde{t}_0]$. Writing
$$[\tilde{t}_0]^{2,4}=t'\otimes \beta,$$
where $\beta$ is the positive generator of  $H^4(\mathbb{P}^3,\mathbb{Z})$, we have now
$$(Id,\phi)_*([\tilde{t}_0]^{2,4})=t'\otimes\phi_*\beta.$$
But since $\phi$ is given by polynomials of degree $p$, the class
$\phi_*\beta$ is  equal to $p$ times the  positive generator
of $H^6(\mathbb{P}^4,\mathbb{Z})$. Hence
$(Id,\phi)_*([\tilde{t}_0]^{2,4})$ is divisible by $p$, and
since $p$ is coprime to $12$, so is $t\otimes j_*\alpha$.

\cqfd

\section{Application to torsion cycles \label{sec3}}
We apply in this section  Koll\'ar's degeneration
 argument to construct interesting torsion algebraic cycles
on smooth projective complex varieties.

In \cite{totaro}, Totaro constructed
examples of cycles  which are homologous to $0$ and
annihilated by the Abel-Jacobi map, but are not algebraically equivalent to
$0$. In fact they are not divisible, while it is well-known
that the groups $CH^k(X)_{alg}$ of
algebraically
equivalent to $0$ cycles  are divisible.
Totaro plays on the factorization
of the map
$$CH^k(X)\rightarrow ( MU^*(X)\underset{MU^*}\otimes{\mathbb Z})^{2k}$$
of (\ref{mapnom})
through algebraic equivalence. Denoting by
$$Griff^k(X)={\mathcal Z}^k(X)_{hom}/{\mathcal Z}^k(X)_{alg}$$ the quotient
of the group of cycles homologous to $0$ by its subgroup
${\mathcal Z}^k(X)_{alg}$, this provides an invariant
$$Griff^k(X)\rightarrow (MU^*(X)\underset{MU^*}\otimes {\mathbb Z})^{2k}$$
with value in the kernel of the map
$\phi^k:(MU^*(X)\underset{MU^*}\otimes{\mathbb Z})^{2k}
\rightarrow H^{2k}(X,{\mathbb Z})$
of (\ref{formulephik}).

These are the invariants used by Totaro to detect non trivial torsion
elements in the Griffiths goup.  The cycles constructed there turn
out to be also annihilated by the Abel-Jacobi map.

The construction we give here will provide torsion cycles annihilated by
the Deligne cycle class (i.e. cohomology class and Abel-Jacobi map)
and also by the Totaro invariants,  but which  are
non divisible, hence non trivial modulo algebraic
equivalence.  A variant of this  construction
(Theorem \ref{th5}) also shows that one can construct such a
cycle in any level of the Hiroshi Saito  filtration on the groups
$CH^*$ (cf \cite{hsaito}).

We shall also construct non trivial
 torsion cycles, which are algebraically equivalent
to $0$ and  annihilated by the
Deligne cycle class map.

The interesting point is that in all three cases, the cycles cannot be detected by
any locally constant  invariant associated to a torsion cycle
$Z\in CH^*(X)_{tors}$
in a smooth complex projective variety $X$. By this, we mean
an invariant which takes value in a locally constant
sheaf on the parameter space for $X$. Indeed,   our
construction provides  torsion  cycles on  smooth
varieties which have parameters,
and we shall see that for some special value of the parameter, our cycles become rationally
equivalent to $0$. Hence any locally constant
invariant attached to them is $0$. In the case of cycles modulo algebraic equivalence,
this is to our knowledge the first
example
of this phenomenon.  Recall that, in contrast, for fixed $X$, torsion
cycles are known to be discrete, by the following lemma (cf \cite{roitman}):

\begin{lemm} (Roitman) If $X$ is projective, $W$ is smooth and connected,
and $\Gamma\in CH^k(W\times
X)$
satisfies the property that for
any $w\in W$, the cycle $\Gamma_*w\in CH^k(X)$ is of torsion,
then $\Gamma_*w$ is in fact constant. Hence, if it vanishes at some point,
it vanishes for any $w$.
\end{lemm}

Our examples show that the last statement is not true if $X$ is allowed
to deform, even staying smooth.

Note finally that, while Bloch-Esnault's \cite{blochesnault}
and Totaro's non divisible  cycles are defined over number fields,
ours might well not be, since we have to restrict to the general point
of a parameter space, which might exclude all the points defined over some
number field.

The general idea of the construction is as follows :
$X$ will be again a hypersurface of degree $p^3$
 in ${\mathbb P}^4$, but we assume that $X$ contains a curve
$C$ of degree
$p$. We arrange things in such a way that the curve $C$
is not divisible by $p$ in $CH^2(X)$, as in the previous section.
This suggests that the cycle
$$pr_1^*c\cdot pr_2^* C$$
on $S\times X$, where the notations are as in section
{\ref{sec2}), is non trivial, which is indeed what we prove.

Let us now give the detailed construction : Fix a prime integer
$p\geq5$. We will consider the following morphism
$$\phi:{\mathbb P}^3\rightarrow{\mathbb P}^4,$$
\begin{eqnarray}\label{phi}\phi(x_0,\ldots,x_3)=(x_0^p,\ldots,x_3^p,f(x_0,\ldots,x_3)),
\end{eqnarray}
where $f$ is a degree $p$ generic homogeneous polynomial. So
$\phi({\mathbb P}^3)=:X_0$ is a hypersurface of ${\mathbb P}^4$ of
degree $p^3$. For any line $l\subset{\mathbb P}^3$, the image
$\phi(l)$ is a curve $C_0$  of degree $p$, and if the line is
conveniently chosen, this curve is smooth, hence is contained in a
smooth hypersurface $X_\infty$  of degree $p^3$.

Consider the pencil $(X_t)_{t\in{\mathbb P}^1}$ generated by
$X_0$ and $X_\infty$. For any $t\in {\mathbb P}^1$, let
$C_t\subset X_t$ be the curve $C_0$, viewed as a curve in $X_t$.
Note that for smooth
$X_t$,  the cohomology class of $C_t$ is divisible
by $p$ : it is equal to $p\alpha$ where $\alpha$ was introduced
in the previous section.
It follows that the Deligne cohomology invariant
$[C_t]_D$ of $C_t$ is also divisible by $p$, since there is the exact
sequence
$$0\rightarrow J^3(X_t)\rightarrow H^4_D(X_t,{\mathbb Z}(2))
\rightarrow H^4(X_t,{\mathbb Z})\rightarrow0,$$
and the intermediate Jacobian $J^3(X_t)$ is  a divisible group.

Let now  $W$ be a smooth projective variety and ${\mathcal L}$
a $p$-torsion  line bundle on $W$. In applications, $W$ will be
either a surface as in the previous section, and
${\mathcal L}$ will not be topologically trivial, with $c_1(\mathcal{L})$ not divisible by $p$ in
$H^2(S,\mathbb{Z})$,
or a curve of genus $>0$, and ${\mathcal L}$
will be topologically trivial, but not trivial.
Let us denote  $c=cl({\mathcal L})\in CH^1(W)$.
Since ${\mathcal L}$ is of $p$-torsion, so is the
Deligne cohomology invariant
$[c]_D\in H^2_D(W,{\mathbb Z}(1))$.

The cycle we will consider is
$$c\times C_t:=pr_1^*c\cdot pr_2^*C_t\in CH^3(W\times X_t).$$
Note that this cycle has vanishing Deligne invariant, since
$$[c\times C_t]_D=pr_1^*[c]_D\cdot_D pr_2^*[C_t]_D,$$
where $\cdot_D$ is the product in Deligne cohomology, and
the left factor is $p$-torsion while the right factor is divisible
by $p$.
This cycle is of $p$-torsion, since $c$
is, and it is algebraically equivalent to $0$ when $c$ is,
for example in the case where $W$ is a curve.

We have now the following :

\begin{theo}\label{theo}1)  For a general  point $t\in{\mathbb P}^1$, the $p$-torsion
cycle
$c\times C_t\in CH^3(W\times X_t)$ is non $0$.

2)  If $c$ is not divisible by $p$ in $CH^1(W)$,  so is
$c\times C_t$ in $CH^3(W\times X_t)$ for general $t$, hence in particular
$c\times C_t$ is non trivial in $Griff^3(W\times X_t)$.
\end{theo}
\noindent {\bf Proof.}
 Let $${\mathcal X}:=
\{(x,t)\in{\mathbb P}^4\times{\mathbb P}^1/x\in X_t\}.$$
${\mathcal X}$ contains the $2$-cycle
 $\Gamma:=C_0\times{\mathbb P}^1$,
and
$W\times {\mathcal X}$ contains the cycle
$$c\times\Gamma:=pr_1^*c\cdot pr_2^*\Gamma.$$
Each $W\times X_t\stackrel{j_t}{\hookrightarrow} W\times{\mathcal X}$ is
the fiber of the natural composite map
$$W\times{\mathcal X}\stackrel{pr_2}{\rightarrow}{\mathcal X}
\stackrel{f}{\rightarrow}{\mathbb P}^1,$$
where $f$ is the map to ${\mathbb P}^1$  given by the pencil.
The restrictions $j_t^*:CH_l(W\times{\mathcal X})
\rightarrow CH_{l-1}(W\times X_t)$ are well defined even if
$W\times{\mathcal X}$ is singular, because
$j_t$ is the inclusion of a Cartier divisor.
We have  $j_t^*(c\times\Gamma)=c\times C_t$ in $CH_i(W\times X_t)$, where
$i=
dim\,W$.
\begin{lemm}\label{lemm} 1) Assume that for a general
complex point $t\in {\mathbb P}^1$,
the restriction $j_t^*(c\times \Gamma)$ is equal to $0$ in
$CH_i(W\times X_t)$. Then for any $t\in {\mathbb P}^1$,
the restriction $j_t^*(c\times \Gamma)$ is equal to $0$ in
$CH_i(W\times X_t)$.

2) Similarly, if for a general complex point $t\in {\mathbb P}^1$,
the restriction $j_t^*(c\times \Gamma)$ is equal to $0$ in
$CH_i(W\times X_t)\otimes \mathbb{Z}/p\mathbb{Z}$, then the same is true for any
$t\in{\mathbb P}^1$.
\end{lemm}

\noindent {\bf Proof.}
 Indeed, the assumption in 1) implies that there exist
a smooth projective curve $D$ and a finite morphism
$r:D\rightarrow\mathbb{P}^1$, such that
denoting by
$$\mathcal{X}_D\stackrel{\tilde{f}}{\rightarrow} D,\,
\mathcal{X}_D\stackrel{\tilde{r}}{\rightarrow}\mathcal{X}
$$ the fibered product
$\mathcal{X}\times_{\mathbb{P}^1}D$, the cycle
$$(Id, \tilde{r})^*(c\times\Gamma)\in CH_{i+1}(W\times\mathcal{X}_D)$$
vanishes on some dense open set $W\times\mathcal{X}_U$, where
$\mathcal{X}_U:=\tilde{f}^{-1}(U)$ for some dense Zariski open set
$U\subset D$.

Similarly in case 2), the cycle
$$(Id, \tilde{r})^*(c\times\Gamma)\in CH_{i+1}(W\times\mathcal{X}_D)$$
will vanish modulo $p$ on some dense open set $W\times\mathcal{X}_U$.

This implies that the cycle $(Id, \tilde{r})^*(c\times\Gamma)$
is supported on fibers of the map $\tilde{f}\circ pr_2:W\times\mathcal{X}_D\rightarrow
D$, that is
\begin{eqnarray}\label{supporte}(Id, \tilde{r})^*(c\times\Gamma)=
\sum_d\tilde{j}_{d*}\gamma_d\,\,\,
{\rm in}\,\,\,CH_{i+1}(W\times\mathcal{X}_D),
\end{eqnarray}
where $\tilde{j}_d$ is the inclusion
of $W\times \tilde{f}^{-1}(d)\cong W\times f^{-1}(r(d))$ in
$W\times\mathcal{X}_D$, and the sum on the right is finite.
Similarly in case 2), we will get the same
equation as in (\ref{supporte}), but modulo $p$.

For any $t\in \mathbb{P}^1$, let $t'\in D$ be such that $r(t')=t$.
Then we have
\begin{eqnarray}\label{eqfin}\tilde{j}_{t'}^*((Id, \tilde{r})^*(c\times\Gamma))=
j_t^*(c\times\Gamma)\,\,{\rm in}\,\,CH_i(W\times X_t).
\end{eqnarray}

But for any $d\in D$, the cycle
$$\tilde{j}_{t'}^*(\tilde{j}_{d*}(\gamma_d))$$
vanishes in $CH_i(W\times X_t)$ since either $t'\not=d$ and the two fibers
do not meet, or $t'=d$, and then $\tilde{j}_{t'}^*\circ\tilde{j}_{d*}$
is by definition the intersection with the class of the Cartier
divisor $ cl({\mathcal O}_{W\times \mathcal{X}_D}(W\times X_{t'})_{\mid W\times X_{t'}})$,
which is trivial.

Combining this with  (\ref{supporte}) and (\ref{eqfin}), we have shown that
$j_t^*(c\times\Gamma)=0\,\,{\rm in}\,\,CH_i(W\times X_t)$ in case 1)
and $j_t^*(c\times\Gamma)=0\,\,{\rm in}\,\,CH_i(W\times X_t)\otimes\mathbb{Z}/p\mathbb{Z}$
in case 2).

\cqfd
We now conclude the proof of Theorem \ref{theo} by contradiction. Assuming
the conclusion of 1) or 2) in the theorem is wrong, we conclude using Lemma \ref{lemm}
that for $t=0$, the cycle $c\times C_0\in CH_i(W\times X)$ is trivial
(resp. is trivial mod. $p$ in case 2) ). This means that we can write
\begin{eqnarray}\label{eq1}
c\times C_0=\sum_i\tau_{i*} div\,\phi_i\,\,{\rm in}\,\,\mathcal{Z}_i(W\times X_0),
\end{eqnarray}
where $W_i$ is normal irreducible of dimension $i+1$, $\phi_i$ is a
non zero  rational function on $W_i$, and $\tau_i:W_i\rightarrow W\times X$
is proper. In case 2), this equality will be true in
$\mathcal{Z}_i(W\times X_0)$ mod $p$.

We observe now that the map $\phi:\mathbb{P}^3\rightarrow X_0$ has the following property :
\begin{lemm}\label{petitlemme}
 For any irreducible closed algebraic subset  $Z\subset \mathbb{P}^3$,
the restriction $\phi_{\mid Z}$ has generic  degree $1$ onto its image.
\end{lemm}
{\bf Proof.} Indeed, if $Z$ is as above, and $\phi$ has degree
$\geq2$ on $Z$, then looking at the first $4$ coordinates
$(x_0^p,\ldots,x_3^p)$ of $\phi$, one concludes that $Z$ has to be
invariant under a group $\mathbb{Z}/p\mathbb{Z}$ acting non
trivially on $Z$, with generator $g$ acting on coordinates by
$$g^*(X_i)=\lambda_iX_i,$$
with $\lambda_i^p=1$. Furthermore, we must have, looking at the last
coordinate of $\phi$: $g^*f=f$ on $Z$. Thus $Z$ has to be contained
in the zero locus of
$$f-g^*f,\,f-(g^2)^*f,\,\ldots,f-(g^{p-1})^*f.$$
But as $f$ is generic, and $p\geq5$, this zero locus consists of the
fixed points of $g$, for any such $g$. Thus $g$ acts trivially on
$Z$, which is a contradiction. \cqfd
 It follows that the
$\tau_i:W_i\rightarrow W\times X_0$ lift to
$\tilde{\tau}_i:W_i\rightarrow W\times{\mathbb{P}^3}$, so that
equation (\ref{eq1}) provides
\begin{eqnarray}\label{eq2}
c\times l=\sum_i\tilde{\tau}_{i*} div\,\phi_i+z\,\,{\rm in}\,\,\mathcal{Z}_i(W\times
\mathbb{P}^3),
\end{eqnarray}
where we recall that $l$ is the line such that $\phi(l)=C_0$ and where $z$
 is a cycle which satisfies :
 \begin{eqnarray}\label{van}(Id,\phi)_*z=0\,\,{\rm in}\,\,\mathcal{Z}_i(W\times X_0).
\end{eqnarray}
 In case 2), this equation becomes :
\begin{eqnarray}\label{eq3}
c\times l=\sum_i\tilde{\tau}_{i*} div\,\phi_i+z+pz'\,\,{\rm in}\,\,\mathcal{Z}_i(W\times
\mathbb{P}^3),
\end{eqnarray}
for some cycle $z'$, where $z$ satisfies property (\ref{van}).
Let $H=cl(\mathcal{O}_{{\mathbb P}^3}(1))\in CH^1(\mathbb{P}^3)$.
We have a map
$$q:CH_i(W\times \mathbb{P}^3)\rightarrow CH_{i-1}(W)=CH^1(W),$$
defined by
$$q(\gamma)=pr_{1*}(\gamma\cdot pr_2^*H).$$
We observe that, because $<l,H>=1$, we have
$$q(c\times l)=c.$$
Applying $q$ to the right-hand sides in equations (\ref{eq2}) and
(\ref{eq3}) gives now
\begin{eqnarray}\label{eq4}
c=q(z)\,\,{\rm in}\,\,CH^1(W),\,\,{\rm resp.}\,\,
c=q(z)\,\,{\rm in}\,\,CH^1(W)\otimes\mathbb{Z}/p\mathbb{Z}\,\,{\rm in \,\,case\,\,2)},
\end{eqnarray}
where $z$ satisfies property (\ref{van}).
We have now the following lemma :
\begin{lemm}\label{easy} If $z\in\mathcal{Z}_i(W\times\mathbb{P}^3)$ satisfies
property (\ref{van}), then $z=0$ in $CH_i(W\times \mathbb{P}^3)$.

\end{lemm}
\noindent {\bf Proof.} The kernel of the map $(Id,\phi)_*$ defined
on cycles is generated by cycles $z_1-z_2$, where $z_i$'s are
effective and $(Id,\phi)_*(z_1)=(Id,\phi)_*(z_2)$. By Lemma
\ref{petitlemme}, we can even assume that the $z_i$ are irreducible
and  $(Id,\phi)(z_1)=(Id,\phi)(z_2)$; looking at the form of the map
$\phi$ given in (\ref{phi}), we see that this is equivalent to the
following :

{\it There exists $\lambda_\cdot=(\lambda_0,\ldots,\lambda_3)$,
where $\lambda_i$ are $p$-th roots of unity, acting  on
$\mathbb{P}^3$  by
$\lambda_\cdot(x_0,\ldots,x_3)=(\lambda_0x_0,\ldots,\lambda_3x_3)$,
such that on $z_1$ the equality
$f(\lambda_\cdot(x_0,\ldots,x_3))=f(x_0,\ldots,x_3)$ is satisfied,
and we have
$$z_2=(Id,\lambda_\cdot)(z_1).$$
}

But the map $(Id,\lambda_\cdot)$ acts as the identity on $CH(W\times\mathbb{P}^3)$,
because its graph is rationally equivalent to the graph of the identity.
Hence $z_2$ is rationally equivalent to $z_1$.

\cqfd
It follows from this lemma that equation (\ref{eq4}) becomes in fact :

\begin{eqnarray}\label{eq5}
c=0\,\,{\rm in}\,\,CH^1(W),\,\,{\rm resp.}\,\,
c=0\,\,{\rm in}\,\,CH^1(W)\otimes\mathbb{Z}/p\mathbb{Z}\,\,{\rm in \,\,case\,\,2)},
\end{eqnarray}
which contradicts the fact that $c\not=0$, resp. $c\not=0$ mod. $p$ in case 2).
This concludes the proof of the theorem.
\cqfd
Theorem \ref{theo}, applied to the case where $\mathcal{L}$ is topologically trivial,
hence $c$ is algebraically equivalent to $0$, provides non trivial torsion cycles
algebraically equivalent to $0$ and annihilated by the Abel-Jacobi map.
Since cycles algebraically equivalent to $0$
are images via correspondences of $0$-cycles homologous to
$0$, this is to be put in contrast with Roitman's theorem \cite{roitman}, which says that the
Albanese map is injective on the torsion part of the groups $CH_0$.

As explained below, these cycles can in fact be specialized to $0$ on a smooth fiber of a flat family
of varieties, so that more generally any invariant of locally constant type has to vanish on them.
Note that Chad Schoen \cite{schoen1} also constructed independently non zero torsion cycles,
algebraically equivalent to $0$, with trivial specialization.

Next, in the case where $W=S$ is a surface and $\mathcal{L}$ is as in section \ref{sec2},
so that $c_1(\mathcal{L})$ is a non divisible class in $H^2(S,\mathbb{Z})$,
Theorem \ref{theo} shows that for general $t$, the cycle $c\times C_t$
is a torsion cycle which is not algebraically equivalent to $0$.

On the other hand,  we can arrange
things so that for some flat deformation of the pair $(C_t,X_t)$
to a pair $(C'_t,X'_t)$, with $X'_t$ smooth,
$c\times C'_t$ becomes rationally equivalent to $0$ on $W\times X'_t$.

Indeed,
in the above construction, we may assume that the initial curve
$C_0$ is  a (singular) plane curve of degree $p$. Such a plane curve
deforms in a flat way to a multiple line
$Z$ of multiplicity $p$ in a plane, and one can then  construct a deformation
of the pair $(C_t,X_t)$ to a pair $(Z,X'_t)$.
One verifies that it is possible to do so with a smooth $X'_t$.
The cycle $z$ associated to $Z$ is divisible by $p$, and it follows
that $c\times z=0$ in $CH^3(W\times X'_t)$.

It follows in particular from this  that the cycle $c\times C_t$ has vanishing
associated Totaro's invariant, since these are locally constant under deformation of
the pair ($X$, torsion-cycle on $X$), and more generally any locally constant invariant
of a torsion cycle
must vanish on it.

To conclude this section, let us note that the same proof can be
used to construct $p$-torsion cycles which are not trivial modulo
algebraic equivalence, (in fact non divisible) and which are in the
$k$-th level of the Hiroshi Saito filtration $F^l_{HSaito}$ on Chow
groups,
 where $k$ can be taken arbitrarily large. Recall that
this filtration is smaller than any existing Bloch-Beilinson filtration, and
is defined as follows :
$F^l_{HSaito}CH^m(X)$ is the subgroup of $CH^m(X)$ which is generated via
correspondences by products $z_1\cdot \ldots\cdot z_l$  of
$l$ cycles homologous to $0$.

Namely, with the same notations as above, consider
the case where $W=S$ and $c$ is not divisible by $p$. Consider the cycle
$$c^k\times C_t^k=(c\times C_t)^k=\Pi_ipr_i^*(c\times C_t)\in CH^{3k}(S^k\times X^k_t).$$
It is a $p$-torsion cycle which lies in
$F^k_{HSaito}CH^{3k}(S^k\times X^k)$, since each
 factor $c\times C_t$ is homologous to
$0$.
Now we have
\begin{theo}\label{th5} For general $t$, the cycle
$c^k\times C_t^k$ is non $0$ in $CH^{3k}(S^k\times X^k)\otimes\mathbb{Z}/p\mathbb{Z}$, hence in particular
it is not algebraically equivalent to $0$.
\end{theo}
{\bf Proof.} We introduce as before the whole family
$$\mathcal{X}^{k}_{/\mathbb{P}^1}:=
\mathcal{X}\times_{\mathbb{P}^1}\ldots\times_{\mathbb{P}^1}\mathcal{X}.$$
On $S^k\times \mathcal{X}^{k}_{/\mathbb{P}^1}$, we have the cycle
$$c^k\times \Gamma_k,$$
where $\Gamma$ was introduced at the beginning of the proof of Theorem \ref{theo}, and
$\Gamma_k:=\Pi_lpr_i^*\Gamma\in CH^{2k}(\mathcal{X}^{k}_{/\mathbb{P}^1})$.

The cycle $c^k\times \Gamma_k$ restricts to $c^k\times C_t^k$ on $S^k\times X_t^k$,
and hence, applying Lemma \ref{lemm}, 2), we conclude that it suffices to show that
the restriction of $c^k\times \Gamma_k$ to the fiber $S^k\times X_0^k$ over $0$
is non trivial modulo $p$.

But this restriction is equal to
$c^k\times C_0^k$, where
$C_0$ is the image of a line $l$ in $\mathbb{P}^3$ via the map $\phi$.
We have now an analogue of Lemma \ref{easy} which works
for the map
$\phi^k:(\mathbb{P}^3)^k\rightarrow X_0^k$, and allows to conclude that
it suffices to show that the cycle $c^k\times l^k$ is not divisible by $p$
in $CH^{3k}(S^k\times (\mathbb{P}^3)^k)$.
But
 the cycle $c^k$ is  not divisible by $p$ in $S^k$, because its
cohomology class is not divisible by $p$, and it follows
by applying the map
$$ CH^{3k}(S^k\times(\mathbb{P}^3)^k)\rightarrow CH^k(S^k),$$
$$\gamma\mapsto pr_{1*}(\gamma\cdot pr_2^*H^{\boxtimes k}),$$
that $c^k\times l^k$ is not divisible by $p$
in $CH^{3k}(S^k\times (\mathbb{P}^3)^k)$.
\cqfd
One may wonder whether a similar construction might allow to construct
examples of
non zero torsion cycles in $F^{n+1}_{HSaito}CH^n(X)$ for some
$n$ and
 some smooth complex projective variety $X$. (Note that it is conjectured by Bloch and Beilinson
that $F^{n+1}_{Hsaito}CH^n(X)_\mathbb{Q}=0$.)

\section{Consequences of the Bloch-Kato conjecture}\label{sec5}

When $F$ is a field and $n \geq 0$ an integer, we denote by $K_n^M (F)$ the $n$-th Milnor $K$-group of $F$.
If a prime $p$ is invertible in $F$ there is a symbol map
$$
K_n^M (F) / p \to H^n (F , \mu_p^{\otimes n})
$$
from the Milnor $K$-theory of $F$ modulo $p$ to the Galois cohomology of $F$ with coefficients
 in the $n$-th power of the Galois group of $p$-th roots of unity. We shall say
 that the {\it Bloch-Kato conjecture is true at the prime $p$}
  when this symbol map is an isomorphism for every $n \geq 0$ and any $F$ of
   characteristic different from $p$
   \cite{BK}. Voevodsky proved that the Bloch-Kato conjecture is true at the prime $2$ \cite{V}.
   He and Rost
   are close to the proof of the Bloch-Kato conjecture at every prime.

Bloch noticed the following striking consequence of the Bloch-Kato conjecture.
Let $X$ be a smooth quasi-projective complex manifold, and $\alpha \in H^n (X,{\mathbb Z})$
 an integral cohomology class. We say that {\it $\alpha$ is supported in codimension $q$} when there exists
  a Zariski closed subset $Y \subset X$ of codimension $q$ such that the restriction of $\alpha$
   to $X-Y$ vanishes.

\begin{theo} \label{th6} (Bloch \cite{B}, end of Lecture 5)  Assume that the Bloch-Kato
 conjecture is true at the prime $p$ and that $\alpha \in H^n (X,{\mathbb Z})$
 is killed by $p$. Then $\alpha$ is supported in codimension one.
 \end{theo}

The Bloch-Kato conjecture gives also some information on the Atiyah-Hirze\-bruch
spectral sequence computing the topological $K$-theory (this fact was first noticed for function fields
 by Thomason).

\begin{theo}\label{th7} Under the same assumption as in Theorem \ref{th6}, the image
of $\alpha$ by every differential $d^r$, $r \geq 2$, is supported in codimension two.
\end{theo}

\noindent {\bf Proof.} We consider the Atiyah-Hirzebruch spectral sequence with ${\mathbb Z} / p$
coefficients:
$$
E_2^{st} (X ; {\mathbb Z} / p) = \left\{ \begin{matrix}
 H^s (X,{\mathbb Z}/p) &\hbox{when $t$ is even} \\ 0 &\hbox{when $t$ is odd,} \end{matrix} \right.
$$
which converges to $K_{\rm top}^{s+t} (X , {\mathbb Z} / p)$. If the Bloch-Kato
 conjecture is true at $p$, it is known \cite{SV} that
  the algebraic $K$-theory of $X$ with ${\mathbb Z} / p$
  coefficients is the abutment of a spectral sequence
$$
'E_2^{ st} (X , {\mathbb Z} / p) = \left\{ \begin{matrix}
 H_{\rm Zar}^s (X,\tau_{\leq i} \,
  R \varepsilon_* \, \mu_p^{\otimes i}) &\hbox{if $t = -2i \leq 0$} \\ 0 &\hbox{otherwise,} \end{matrix}
   \right.
$$
where $\varepsilon : X_{\hbox{\footnotesize \'et}} \to X_{\rm Zar}$ is
the natural morphism from the big \'etale site of $X$ to its
 big Zariski site, and $\tau_{\leq i}$ is the good truncation. There is a morphism of spectral sequences
$$
'E_r^{ st} (X,{\mathbb Z}/p) \to E_r^{ st} (X,{\mathbb Z}/p) \, .
$$
Consider the differential
$$
d^r : E_r^{n,-2n} (X,{\mathbb Z}/p) \to E_r^{n+r,-2n-r+1} (X,{\mathbb Z}/p) \, .
$$
The group
$$
'E_2^{n,-2n} (X,{\mathbb Z}/p) = H_{\rm Zar}^n (X , \tau_{\leq n} \, R \varepsilon_* \, \mu_p^{\otimes n})
= H_{\rm Zar}^n (X , R \varepsilon_* \, \mu_p^{\otimes n}) =
 H_{\hbox{\footnotesize \'et}}^n (X , \mu_p^{\otimes n})
$$
maps isomorphically to
$$
E_r^{n,-2n} (X,{\mathbb Z}/p) = H^n (X,{\mathbb Z}/p) \, .
$$
On the other hand, if $i = n + \frac{r-1}{2}$, the spectral sequence
 computing the hypercohomology of $\tau_{\leq i} \, R \varepsilon_* \, \mu_p^{\otimes i}$ is
$$
''E_2^{st} =  \left\{ \begin{matrix} H_{\rm Zar}^s (X, R^t \, \varepsilon_* \, \mu_p^{\otimes i})
 &\hbox{when $0 \leq t  \leq i$} \\ 0 &\hbox{otherwise,} \end{matrix} \right.
$$
with abutment $H_{\rm Zar}^{s+t} (X , \tau_{\leq i} \, R \varepsilon_* \, \mu_p^{\otimes i})$.
When $s+t = n+r$ we find $''E_2^{ st} = 0$ unless $n+r-s \leq i = n + \frac{r-1}{2}$, {\it i.e.}
$$
r \leq 2 \, s - 1 \, .
$$
According to Bloch-Ogus \cite{BO}, any class in $''E_2^{ st}$ is supported in codimension $s$.
Therefore, any class in $H_{\rm Zar}^{n+r} (X , \tau_{\leq i} \, R \varepsilon_* \, \mu_p^{\otimes i})$
is supported in codimension $\frac{r+1}{2} \geq 2$.

Therefore, given $x \in H^n (X,{\mathbb Z} / p) \simeq \, 'E_2^{n,-2n} (X , {\mathbb Z} / p)$, we
can find a Zariski closed subset $Y \subset X$, ${\rm codim} \, (Y) \geq 2$, such that,
 for every $r \geq 2$, the restriction of $'d^r (x) \in \, 'E_2^{n+r,-2n-r+1} (X,{\mathbb Z} / p)$
 to $X-Y$ vanishes. Therefore the same is true for $d^r (x)$.

If $\alpha \in H^n (X,{\mathbb Z})$ is such that $p \alpha = 0$, it lies in the image of
the Bockstein homomorphism
$$
\beta : H^{n-1} (X,{\mathbb Z} / p) \to H^n (X , {\mathbb Z}) \, .
$$
Therefore $\alpha = \beta (x)$, and $d^r (\alpha) = \beta (d^r (x))$ is supported in codimension two.

C.Soul\'e: CNRS and IHES, 35 Route de Chartres, 91440 Bures sur Yvette.  soule@ihes.fr

\medskip

C.Voisin: CNRS, UMR 7586, Institut de Math\'{e}matiques de Jussieu, 175 rue du Chevaleret, 75013 Paris.
voisin@math.jussieu.fr

\end{document}